\documentclass[12pt,a4paper]{article}%
\usepackage{amsmath, amsfonts, amsthm,color, latexsym, amssymb}
\usepackage{amscd}
\usepackage{amsmath}
\usepackage{amsfonts}
\usepackage{amssymb}
\usepackage{graphicx}%
\providecommand{\U}[1]{\protect\rule{.1in}{.1in}}
\newtheorem{theorem}{Theorem}[section]

\newtheorem{lemma}[theorem]{Lemma}

\begin{document}

\title{\textsc{Dominated polynomials on infinite dimensional spaces}}
\author{Geraldo Botelho\thanks{Supported by CNPq Grant 620108/2008-8.}, Daniel Pellegrino\thanks{Supported by by
INCT-Matem\'{a}tica, CNPq Grants 620108/2008-8, 308084/2006-3 and
471686/2008-5.}~ and Pilar Rueda\thanks{Supported by Ministerio de
Ciencia e Innovaci\'{o}n MTM2008-03211/MTM.\hfill\newline2000
}}
\date{}
\maketitle

\begin{abstract}
The aim of this paper is to prove a stronger version of a conjecture
posed in \cite{PAMS} on the existence of non-dominated scalar-valued
$m$-homogeneous polynomials, $m \geq 3$, on arbitrary infinite
dimensional Banach spaces.

\end{abstract}

\vspace*{-1.0em}

\section*{Introduction}

The theory of absolutely summing multilinear mappings and homogeneous polynomials between Banach spaces, which was first outlined by Pietsch \cite{pietsch}, studies multilinear and polynomial generalizations of the very successful theory of absolutely summing linear operators. The theory has been developed by several authors and among the advances obtained thus far we mention: Pietsch-type domination/factorization theorems (\cite{jfa1}, Geiss \cite{geiss}, P\'erez-Garc\'ia \cite{davidarchiv}), different types of absolutely summing multilinear mappings and polynomials (Achour and Mezrag \cite{achour}, Carando and Dimant \cite{cd}, \c Caliskan and Pellegrino \cite{erhan}, Dimant \cite{dimant}, Pellegrino and Souza \cite{marcela}, P\'erez-Garc\'ia \cite{davidarchiv}), Grothendieck-type theorems (Bombal et al \cite{QJM}, P\'erez-Garc\'ia and Villanueva \cite{jmaa}), coincidence/inclusion/composition theorems (Alencar and Matos \cite{am}, Botelho et al \cite{botpams}, P\'erez-Garc\'ia \cite{monat, davidstudia}, Popa \cite{popa}), connections with the geometry of Banach spaces (\cite{irish}, Floret and Matos \cite{Floret}, Mel\'endez and Tonge \cite{MT}, \cite{danielstudia}, P\'erez-Garc\'ia \cite{fenn}), interplay with other multi-ideals and polynomial ideals (Botelho et al \cite{botstudia}, \cite{PRIMS}, Cilia and Guti\'errez \cite{cg}, Jarchow et al \cite{Jar}, Matos \cite{nuclear}), estimates for absolutely summing norms (Aron et al \cite{alrt}, \cite{indagnovo}, Choi et al \cite{choi}, Defant and Sevilla-Peris \cite{sevilla}, Zalduendo \cite{zalduendo}), extensions of the theory to more general nonlinear mappings (Junek, Matos and Pellegrino \cite{junek}, Matos \cite{anais, collec}, Matos and Pellegrino \cite{matosdaniel}). \\
\indent Let us sketch some of the applications this theory has produced. The following results were obtained with the help of the theory of absolutely summing multilinear mappings: (i) for every $n \in \mathbb{N}$, a tensor norm of order $n$ constructed by P\'erez-Garc\'ia and Villanueva \cite{jmaa} is shown in Defant and P\'erez-Garc\'ia \cite{andreas david} to preserve unconditionality for ${\cal L}_p$-spaces; (ii) Defant et al \cite{defant et al} provides optimal estimates for the width of Bohr's strip for Dirichlet series in infinite-dimensional Banach spaces; (iii) applications to quantum information theory are obtained by P\'erez-Garc\'ia et al \cite{david et al}, where it is proved that, contrary to the bipartite case, tripartite Bell inequalities can be unboundedly violated; (iv) in Jarchow et al \cite{Jar}, the existence of Hahn-Banach-type extension theorems for multilinear forms is strongly connected to structural properties of the underlying spaces.

One of the central notions of this theory is that of dominated homogeneous polynomial. A continuous $m$-homogeneous polynomial $P$ from the Banach space $X$ to the
Banach space $Y$ is {\it $r$-dominated} if $(P(x_{j}))_{j=1}^{\infty}$ is $\frac
{r}{m}$-summable in $Y$ whenever $(x_{j})_{j=1}^{\infty}$ is weakly
$r$-summable in $X$.

 The following conjecture was posed in \cite{PAMS}:

\medskip

\noindent\textbf{Conjecture 1.} There is no infinite dimensional Banach space
$X$ such that for every $m \in\mathbb{N}$ and every $r \geq1$, any continuous
scalar-valued $m$-homogeneous polynomial on $X$ is $r$-dominated.

\medskip

It is known that the conjecture holds true for Banach spaces with
unconditional basis (see \cite[Theorem 3.2]{PAMS}). The status of the problem was changed by the proof of
a stronger result for multilinear forms: in Jarchow et al
\cite[Lemma 5.4]{Jar} (see also \cite[Proposition 3.2]{cot-inf}), it is shown that for every $m\geq3$ and every $r\geq1$, on every infinite dimensional Banach space there exists a continuous $m$-linear form that fails to be $r$-dominated. Although not solving Conjecture 1 (to prove Conjecture 1 we need a {\it symmetric} non-$r$-dominated $m$-linear form), this result indicates that a result stronger than Conjecture 1 should be pursued for polynomials:

\medskip

\noindent\textbf{Conjecture 2.} For every infinite dimensional Banach space
$X$, every $m \geq 3$ and every $r \geq1$, there is a continuous
scalar-valued $m$-homogeneous polynomial on $X$ that fails to be $r$-dominated.

\medskip

In Section \ref{conjecture} we solve Conjecture 2 (hence Conjecture 1) in the positive.

\section{Notation}

Throughout this paper $n$ and $m$ are positive integers, $X$ and $Y$ will
stand for Banach spaces over $\mathbb{K}=\mathbb{R}$ or $\mathbb{C}$. The
Banach space of all continuous $m$-homogeneous polynomials $P\colon
X\longrightarrow Y,$ with the $\sup$ norm, is denoted by $\mathcal{P}%
(^{m}X;Y)$ ($\mathcal{L}(X;Y)$ if $m=1$). When $m=1$ and $Y=\mathbb{K}$ we
write $X^{\ast}$ to denote the topological dual of $X$. The closed unit
ball of $X$ is represented by $B_{X}.$ For details on the theory of polynomials between Banach spaces we
refer to \cite{Di, Mu}.

Given $r\in\lbrack0,\infty)$, let $\ell_{r}(X)$ be the Banach ($r$-Banach if
$0 < r < 1$) space of all absolutely $r$-summable sequences $(x_{j}%
)_{j=1}^{\infty}$ in $X$ with the norm $\Vert(x_{j})_{j=1}^{\infty}\Vert
_{r}=(\sum_{j=1}^{\infty}\Vert x_{j}\Vert^{r})^{1/r}$. We denote by $\ell
_{r}^{w}(X)$ the Banach ($r$-Banach if $0 < r < 1$) space of all weakly
$r$-summable sequences $(x_{j})_{j=1}^{\infty}$ in $X$ with the norm
$\Vert(x_{j})_{j=1}^{\infty}\Vert_{w,\,r}=\sup_{\varphi\,\in\,B_{X^{\ast}}%
}\Vert(\varphi(x_{j}))_{j=1}^{\infty}\Vert_{r}$.

Let $p,q>0.$ An $m$-homogeneous polynomial $P\in\mathcal{P}(^{m}X;Y)$ is
absolutely $(p;q)$-summing if $\left(  P(x_{j})\right)  _{j=1}^{\infty}\in
\ell_{p}(Y)$ whenever $\left(  x_{j}\right)  _{j=1}^{\infty}\in\ell_{q}%
^{w}(X).$ It is well-known that $P$ is absolutely $(p;q)$-summing if and only
if there is a constant $C\geq0$ so that%
\[
\left(
{\displaystyle\sum_{j=1}^{n}}
\left\Vert P(x_{j})\right\Vert ^{p}\right)  ^{\frac{1}{p}}\leq C\left(
\Vert(x_{j})_{j=1}^{n}\Vert_{w,\,q}\right)  ^{m}%
\]
for every $x_{1},\ldots,x_{n}\in X$ and $n$ positive integer. The infimum of
such $C$ is denoted by $\Vert P\Vert_{as(p;q)}.$ The space of all
absolutely $(p;q)$-summing $m$-homogeneous polynomials from $X$ to $Y$ is
denoted by $\mathcal{P}_{as(p;q)}(^{m}X;Y)$ and $\Vert\cdot\Vert_{as(p;q)}$ is
a complete norm ($p$-norm if $p<1$) on $\mathcal{P}_{as(p;q)}(^{m}X;Y)$.

An $m$-homogeneous polynomial $P\in\mathcal{P}(^{m}X;Y)$ is said to
be $r$\textit{-dominated} if it is absolutely
$(\frac{r}{m};r)$-summing. In this case we write
$\mathcal{P}_{d,r}(^{m}X;Y)$ and $\Vert\cdot\Vert_{d,r}$ instead of
$\mathcal{P}_{as(\frac{r}{m};r)}(^{m}X;Y)$ and $\Vert\cdot\Vert
_{as(\frac{r}{m};r)}$. As usual  we write $\mathcal{P}_{d,r}(^{m}X)$
and $\mathcal{P}(^{m}X)$ when $Y=\mathbb{K}$. The definition (and
notation) of $r$-dominated multilinear mappings is analogous (for
the notation just replace $\cal P$ by $\cal L$). The symbol ${\cal
L}^s$ means that only symmetric multilinear mappings are considered.
For details we refer to \cite{PAMS, Jar}.

\section{The proof of Conjecture 2}\label{conjecture}

It is well-known that a homogeneous polynomial is $p$-dominated if and only if
its associated symmetric multilinear map is also $p$-dominated (see
\cite[Theorem 6]{MT}). However, the multilinear and polynomial parts of our problem are not so tightly connected, as is made clear by the following example of a non-$p$-dominated multilinear form whose symmetrization is $p$-dominated: 
let $$T:\ell
_{2}\times\ell_{2}\rightarrow\mathbb{K}~,~
T(x,y)=%
{\displaystyle\sum\limits_{j=1}^{\infty}}
x_{j}y_{j+1}-%
{\displaystyle\sum\limits_{j=1}^{\infty}}
x_{j+1}y_{j}.
$$
Note that $T(e_{j},e_{j+1})=1$ for every $j$. So, $T$ fails to be
$p$-dominated (regardless of the $p\geq1$). On the other hand the
symmetrization of $T$ is zero. So, even knowing that the multilinear
counterpart of the conjecture is true, it is in principle not clear that the same
holds for polynomials.

To prove Conjecture 2 we need the following well-known results:

\begin{lemma}
\label{note} {\rm ~ \\
{\rm (a)} \cite[Propositions 41(b) and 46(a)]{note}} \label{dan}If
$\mathcal{P}_{d,r}(^{n}X;Y)=\mathcal{P}(^{n}X;Y)$ then $\mathcal{P}%
_{d,r}(^{m}X;Y)=\mathcal{P}(^{m}X;Y)$ for every $m\leq n.$\\
{\rm (b) (\cite[Proposition 3.4]{jfa1} and \cite[Theorem 2.8]{DJT})} $\mathcal{P}_{d,r}(^{n}X;Y) \subseteq \mathcal{P}_{d,q}(^{n}X;Y)$ if $r \leq q$.
\end{lemma}

The following result has its proof inspired on the proof of \cite[Lemma
5.4]{Jar}; its final form was kindly suggested by the referee.

\begin{theorem}\label{theorem}
If $X$ is infinite dimensional and $\mathcal{P}_{d,2}(^{2}X) = \mathcal{P}(^{2}X)$,
then $\mathcal{P}_{d,r}(^{m}X)\neq\mathcal{P}(^{m}X)$
for every $m\geq3$ and every $r\geq1.$
\end{theorem}

\begin{proof} By Lemma \ref{note}(b) we may assume $r\geq2$ and by Lemma
\ref{note}(a) we just need to prove the case $m=3$.
Since $X$ is infinite dimensional, by Dvoretzky's theorem there is a sequence $(E_n)_{n=1}^\infty$ of subspaces of $X$ such that $X/E_n$ is 2-isomorphic to $\ell_2^n$ for every $n$ (details can be found in the proof of \cite[Proposition 19.17(b)]{DJT}). So we can consider isomorphisms $j_n \colon X/E_n \longrightarrow \ell_2^n$ such that $\|j_n\| \leq 4$ and $\|j_n^{-1}\| = \frac12$ for every $n$. Letting $\pi_n \colon X \longrightarrow X/E_n$ be the corresponding quotient maps and defining $q_n := j_n \circ \pi_n$ we have that
$$\|q_n \| \leq 4 {\rm ~~and~~} B_{\ell_2^n} \subseteq q_n(B_X)~~{\rm for~every~}n. $$
By assumption we know that
${\mathcal L}^s_{d,2}(^2X)={\mathcal L}^s(^2X)$, so by the open
mapping theorem there exists $C>0$, depending only on $X$, such that
$\|A\|_{d,2}\leq C \|A\|$ for every $A \in {\cal L}^s(^2X)$. Given $n \in \mathbb{N}$, defining
$$A \colon X \times X \longrightarrow \mathbb{K}~,~A(x,y) = \langle q_n(x), q_n(y) \rangle, $$
we have that $A$ is a continuous symmetric (in the complex case we consider its symmetrization) bilinear form on $X$ and $\|A\| \leq 16$ because $\|q_n\| \leq 4$. Given $x_1, \ldots, x_k \in X$,
\begin{eqnarray*}
\sum_{j=1}^k \|q_n(x_j)\|_2^2 &=& \sum_{j=1}^k \langle q_n(x_j), q_n(x_j) \rangle ~=~ \sum_{j=1}^k |A(x_j,x_j)|\\
&\leq& \|A\|_{d,2} \|(x_j)_{j=1}^k\|_{w,2}^2 ~\leq~ C \|A\|
\|(x_j)_{j=1}^k\|_{w,2}^2 \\
&\leq& 16C \|(x_j)_{j=1}^k\|_{w,2}^2,
\end{eqnarray*}
showing that $\|q_n\|_{as(2,2)} \leq 4\sqrt{C}$ for every $n$.

For each positive integer $n$, consider the canonical embedding $i_n \colon \ell_2^n \longrightarrow \ell_\infty^n$. It is plain that $\|i_n\| = 1$ for every $n$ and $\lim_{n\rightarrow\infty}\left\Vert i_{n}\right\Vert _{as(r;r)}=\infty$. Again by Dvoretzky's theorem, there is an $n$-dimensional
 subspace $X_n$ of $X$ and an isomorphism $k_n \colon X_n \longrightarrow \ell_2^n$ such that $\|k_n\| \leq 1$ and $\|k_n^{-1}\| \leq 2$ (for every $n$). Since $\ell_\infty^n$
 is an injective space, there is a norm-preserving extension $u_n$ of $i_n \circ k_n$ to the whole of $X$. Let us see that $\sup_{n}\left\Vert u_{n}\right\Vert _{as(r;r)}=\infty$. Given $M > 0$, since $\sup_n\|i_n\|_{as(r,r)}=\infty$, there exists $n$ such that $\|i_n\|_{as(r,r)}>M$.
Then, there exist $m\in \mathbb{N}$ and $x_1,\ldots,x_m\in \ell_2^n$
such that
\begin{eqnarray}\label{e}
\left(\sum_{j=1}^m\|i_n(x_j)\|_\infty^r\right)^{1/r}>M\|(x_j)_{j=1}^m\|_{w,r}.
\end{eqnarray}
Since the restriction of $u_n$ to $X_n$ is $i_n\circ k_n$, for
any $x\in \ell_2^n$ we have
$$
\|i_n(x)\|_\infty=\|i_n\circ
k_n(k_n^{-1}(x))\|_\infty=\|u_n(k_n^{-1}(x))\|_\infty.
$$
From (\ref{e}) it follows that
\begin{eqnarray*}
\left(\sum_{j=1}^m\|u_n(k_n^{-1}(x_j))\|_\infty^r \right)^{1/r}&>&
M\|(x_j)_{j=1}^m\|_{w,r}\\
&=& \frac{M}{2}\sup_{\varphi \in 2B_{(\ell_2^n)^*}}\left(\sum_{j=1}^m
\left|\varphi(x_j)\right|^r\right)^{1/r}\\
&\stackrel{(*)}{\geq}& \frac{M}{2}\sup_{\psi \in B_{(X_n)^*}}\left(\sum_{j=1}^m
\left|\psi(k_n^{-1}(x_j))\right|^r\right)^{1/r}\\
&=&\frac M2\|(k_n^{-1}(x_j))_{j=1}^m\|_{w,r}.
\end{eqnarray*}
$(*)$ $\psi \in B_{(X_n)^*} \Longrightarrow \|\psi \circ k_n^{-1}\| \leq 2 \Longrightarrow \psi \circ k_n^{-1} \in 2B_{(\ell_2^n)^*}$.\\
Hence, $\sup_n\|u_n\|_{as(r,r)}=\infty$. Following the steps of the
 proof of \cite[Lemma 5.4]{Jar}, $\ell_{\infty}^{n}$ can be identified with the diagonal of $\ell_{2}^{n}{\otimes
}_{\varepsilon}\ell_{2}^{n}$, so that $u_n$ can be regarded as a linear operator from $X$ to the diagonal of $\ell_{2}^{n}{\otimes
}_{\varepsilon}\ell_{2}^{n}$ whose norm does not depend on $n$. From
$$\ell_{2}^{n}{\otimes}_{\varepsilon}\ell_{2}^{n} \stackrel{(**)}{=}
((\ell_{2}^{n})^*{\otimes}_{\pi}(\ell_{2}^{n})^*)^* = (\ell_{2}^{n}{\otimes}_{\pi}\ell_{2}^{n})^* = {\cal L}(^2 \ell_2^n)$$
(for $(**)$ see \cite[Proposition 4.1(1)]{DF}), $u_n$ can now be regarded as a linear operator from $X$ to ${\cal L}(^2 \ell_2^n)$ such that $u_n(x)$ is a symmetric bilinear form (elements of the diagonal of the tensor product are symmetric tensors, to which correspond symmetric bilinear forms) for every $x \in X$. 

Consider the 3-linear form
\[
T_{n}:X\times X\times X\longrightarrow\mathbb{K}~,~T_{n}(x,y,z)=u_{n}(x)(q_{n}(y),q_{n}(z)),
\]
its symmetrization $T_n^s \in {\cal L}^s(^3X)$,
\[
{T}_{n}^s(x,y,z)=\frac{1}{3}\left[  T_{n}(x,y,z)+T_{n}(y,x,z)+T_{n}%
(z,x,y)\right],
\]
and the linear operator $\widetilde{T}_{n}^s \colon X \longrightarrow {\cal L}(^2X)$ associated to $T_n^s$,
\begin{eqnarray*}
\widetilde{T}_n^s(x)(y,z)\!\!&=&\!\! \frac13\left[T_n(x,y,z) + T_n(y,x,z) + T_n(z,x,y) \right]\\
\!\!&=&\!\! \frac{1}{3}\left[  u_{n}(x)(q_{n}(y),q_{n}(z)) +
u_{n}(y)(q_{n}(x),q_{n}(z)) + u_{n}(z)(q_{n}(x),q_{n}(y))\right].
\end{eqnarray*}
Write $3\widetilde{T}_n^s = A_n + B_n + C_n$, where
$$A_n(x)(y,z) =  u_{n}(x)(q_{n}(y),q_{n}(z)),$$
$$\hspace*{1.9em}B_n(x)(y,z) = u_{n}(y)(q_{n}(x),q_{n}(z)),~{\rm and}$$
$$C_n(x)(y,z) = u_{n}(z)(q_{n}(x),q_{n}(y)).$$
Given $x_1, \ldots, x_k \in X$, 
\begin{eqnarray*}
\sum_{j=1}^k \|B_n(x_j)\|^2&=&\sum_{j=1}^k \sup_{y,z \in B_X}|u_{n}(y)(q_{n}(x_j),q_{n}(z))|^2 \\& \leq & \sum_{j=1}^k \|u_{n}\|^2 \|q_{n}(x_j)\|^2\|q_{n}\|^2 \\
& \leq & \sup_n \|u_n\|^2 \sup_n \|q_{n}\|^2 \|q_{n}\|_{as(2,2)}^2\|(x_j)_{j=1}^k\|_{w,2}^2,
\end{eqnarray*}
proving that $B_n$ is $2$-summing (hence $r$-summing because
$r\geq2$) and its 2-summing norm (hence its $r$-summing norm) is
controlled by $16\sqrt{C}\sup_n \|u_n\|$. It is clear that the same
happens to $C_n$, so
$$\sup_n \|B_n\|_{as(r,r)} < \infty ~~{\rm and}~~\sup_n \|C_n\|_{as(r,r)} < \infty.$$
Let us prove that $\sup_n \|A_n\|_{as(r,r)} = \infty$. Given $M'
> 0$, from $\sup_n \|u_n\|_{as(r,r)} = \infty$ there is $n$ such
that $\|u_n\|_{as(r,r)} > M'$. So, there are $m \in \mathbb{N}$ and
$y_1, \ldots, y_m \in X$ such that
$$\left(\sum_{j=1}^m \|u_n(y_j)\|^r \right)^{\frac{1}{r}} > M' \|(y_j)_{j=1}^m\|_{w,r}. $$ Using that $B_{\ell_2^n} \subseteq q_n(B_X)$, we have
\begin{eqnarray*}
\|u_n(x)\| &=& \sup_{\lambda_1, \lambda_2 \in B_{\ell_2^n}} |u_n(x)(\lambda_1, \lambda_2)|\\
&\leq& \sup_{y,z \in B_X}|u_n(x)(q_n(y), q_n(z))|\\
&=& \sup_{y,z \in B_X} |A_n(x)(y,z)|~=~\|A_n(x)\|,
\end{eqnarray*}
for every $x \in X$. So,
$$\left(\sum_{j=1}^m \|A_n(y_j)\|^r \right)^{\frac{1}{r}} \geq \left(\sum_{j=1}^m \|u_n(y_j)\|^r \right)^{\frac{1}{r}} > M' \|(y_j)_{j=1}^m\|_{w,r}, $$
proving that $\|A_n\|_{as(r,r)} > M'$, hence
$\sup_n \|A_n\|_{as(r,r)} = \infty$. From $3\widetilde{T}_n^s = A_n + B_n + C_n$, it follows that $\sup_n \|\widetilde T_n^s\|_{as(r,r)} = \infty$. Hence $\sup_n \|T_n^s\|_{d,r} = \infty$ because $\|\widetilde{T}_n^s\|_{as(r,r)} \leq \|T_n^s\|_{d,r}$ (see, e.g., \cite[Lemma 3.4]{irish}). But $\sup_n \|T_n^s\| < \infty$ because $\sup_n \|u_n\| < \infty$ and $\sup_n \|q_n\| < \infty$, so the open mapping theorem yields that ${\cal L}^s(^3X) \neq {\cal L}_{d,r}^s(^3X)$. Therefore ${\cal P}(^3X) \neq {\cal P}_{d,r}(^3X)$.
\end{proof}

To complete the proof of Conjecture 2 we need another well-known result (for complex scalars the proof appeared in \cite[Corollary 3.2]{Floret}, the general case can be found in \cite[Proposition 13]{cg}):

\begin{lemma}
\label{novo-dez}$\mathcal{P}_{d,r}(^{2}X)=\mathcal{P}_{d,2}(^{2}X)$ for every Banach space $X$ and $r\geq2$.%
\end{lemma}

\begin{theorem}\label{tttt}
Let $m\geq3,$ $r\geq1$ and $X$ be an infinite dimensional Banach
space. Then%
\begin{equation}
\mathcal{P}_{d,r}(^{m}X)\neq\mathcal{P}(^{m}X)\text{.}
\label{nbvc}%
\end{equation}
Moreover, if $\mathcal{P}_{d,2}(^{2}X)\neq\mathcal{P}(^{2}%
X)$, then {\rm (\ref{nbvc})} holds for every $m\geq2$ and $r\geq1.$
\end{theorem}

\begin{proof} If $\mathcal{P}_{d,2}(^{2}X)= \mathcal{P}(^{2}X)$, Theorem \ref{theorem} gives the result.
If $\mathcal{P}_{d,2}(^{2}X)\neq\mathcal{P}(^{2}X)$, Lemma \ref{novo-dez}
gives that $\mathcal{P}_{d,r}(^{2}X)\neq%
\mathcal{P}(^{2}X)$ for every $r\geq 2$. Hence $\mathcal{P}_{d,r}(^{2}X)\neq\mathcal{P}(^{2}X)$ for every $r\geq 1$ by Lemma \ref{note}(b). Lemma \ref{note}(a) completes the result for every $m\geq2$.
\end{proof}

\medskip

\noindent {\bf Acknowledgements.} The authors thank N. Kalton and the referee for their helpful suggestions.

\vspace*{1em} \noindent[Geraldo Botelho] Faculdade de Matem\'atica,
Universidade Federal de Uberl\^andia, 38.400-902 - Uberl\^andia, Brazil,
e-mail: botelho@ufu.br.

\medskip

\noindent[Daniel Pellegrino] Departamento de Matem\'atica, Universidade
Federal da Para\'iba, 58.051-900 - Jo\~ao Pessoa, Brazil, e-mail: dmpellegrino@gmail.com.

\medskip

\noindent[Pilar Rueda] Departamento de An\'alisis Matem\'atico, Universidad de
Valencia, 46.100 Burjasot - Valencia, Spain, e-mail: pilar.rueda@uv.es.

\end{document}